\documentclass[a4paper,11pt]{article}
\usepackage{pstricks}
\usepackage{amsmath}
\usepackage{amsfonts}
\usepackage{amscd}
\usepackage{amssymb}
\usepackage[all]{xy}

\newcommand{\al}{{\alpha}}

\newcommand{\om}{{\omega}}

\newcommand{\la}{{\lambda}}

\newcommand{\rar}{\rightarrow}

\newcommand{\NI}{{\noindent}}

\newcommand{\bk}{\hfill $\Box$}

\newtheorem{theorem}{Theorem}[section]
\newtheorem{thm}[theorem]{Theorem}
\newtheorem{cor}[theorem]{Corollary}
\newtheorem{defin}[theorem]{Definition}
\newtheorem{rem}[theorem]{Remark}
\newtheorem{lemma}[theorem]{Lemma}
\newtheorem{prop}[theorem]{Proposition}
\newtheorem{ex}[theorem]{Example}

\newtheorem{con}[theorem]{Conjecture}

\numberwithin{equation}{section}

\title{\bf On vector fields having properties of Reeb fields}
\author{Bogus\l\/aw Hajduk  and Rafa\l\ Walczak }

\date{}

\begin{document}

\maketitle

\begin{abstract}

We study constructions of vector fields
with properties which are characteristic to Reeb vector fields
of contact forms. In particular, we prove that all closed oriented
odd-dimensional manifold have geodesible vector fields.

\noindent {\bf Keywords:} Reeb field, presymplectic form, contact form, geodesible field

\noindent {\bf AMS classification(2010)}: Primary 53D15, 53D35; Secondary 57R25

\end{abstract}

\bigskip

\section{Introduction}

One of the major open questions in contact topology can be formulated as follows.

\begin{con}\label{con1} Does any closed, oriented, odd-dimensional smooth manifold $M$
having a non-degenerate two-form admit a contact form?
\end{con}

\noindent Here non-degeneracy means that the rank of the form is maximal at
each point, hence it is $2n$ if $\dim M=2n+1.$ A non-degenerate two
form $\omega$ exists if and only if the structure group of the
tangent bundle $TM$ has a  reduction to the group $U(n)\subset
SO(2n+1).$ The kernel of such a form is a 1-dimensional subbundle
$\mathbf R$ of $TM,$ hence it gives a 1-dimensional foliation $\mathfrak R$
of $M.$ Properties of Reeb fields are studied intensively in contact topology
 and some sophisticated analytical techniques are used for this purpose, as for instance
symplectic field theory (SFT) or contact homology.
In the present note we deal with
properties of Reeb foliations which can be detected by topological
methods. The general question, which is our motivation, is the following:
given a structure on $M$ which arises in the presence of a contact
form on $M,$  find constructive methods to build
such a structure, at least on some classes of manifolds.

In the paper we consider the following two
properties which are satisfied by the Reeb foliation $\mathfrak R$
 of any contact form:

\medskip
\noindent {\bf Property A:} $\mathfrak R$ admits a connection (equivalently,
it is geodesible);

\noindent {\bf Property B.} the basic cohomology class
of exterior derivative of a connection on $\mathfrak R$ is non-zero.

\medskip
In the next section we show that in fact these properties are
valid for Reeb fields (hence for foliations they define) of
contact forms.
Our main result is that on any closed oriented odd-dimensional
manifold $M$ a vector field having both properties can be
constructed starting from an open book decomposition of $M.$

There are obstructions for a foliation to have
Property A found by Dennis Sullivan \cite{S}. However,
Herman Gluck announced \cite{G} that any closed manifold of
odd dimension admits a foliation with Property A, but
according to our knowledge, the
proof was never published. In this paper, we give a simple proof of
this hypothesis, based on existence of open book decompositions,
cf. Theorem \ref{main}.

We give also an example showing that a presymplectic form (i.e.,
non-degenerate and closed) does not need to satisfy neither A nor B.

\section{Properties of Reeb vector fields}

Throughout this paper we assume all manifolds to be smooth, closed
and oriented. Non-zero vector fields are denoted by $V,R,..,$ and
the 1-dimensional foliations generated by the fields
will be denoted by the corresponding Gothic letter
$\mathfrak{V}, \mathfrak{R}...$ The
contraction of a form $\eta$ with a vector field $V$
will be denoted by $\iota_V\eta .$

 A {\it contact form} on $M^{2n+1}$ is a one-form $\alpha$ such
that $\alpha \wedge (d\alpha)^n>0.$ With every contact form we
associate a vector field called
its {\it Reeb field } R. It is uniquely defined by
two conditions $\alpha(R)=1$ and $\iota_Rd\alpha=0.$ The
condition that defines a contact form implies that
$d\alpha$ is non-degenerate (i.e. it has rank $2n$).

A closed 2-form $\omega$ such that $\om$ is
non-degenerate is called {\it presymplectic.} In this case
we have the {\it Reeb foliation} $\mathfrak{R}$
of $\om$ defined by $\iota_{\mathfrak{R}}\om=0.$ The
$2n$-form $\om^n$ defines an orientation on any subbundle
transversal to $\mathfrak{R},$ hence $\mathfrak{R}$ has
non-zero sections (non-vanishing vector fields tangent
to $\mathfrak R).$ All such fields will be called Reeb fields.
 By definition, $\iota_R\om=0.$

By \cite{MD}, any non-degenerate two-form can be deformed to
a presymplectic form. Thus Conjecture \ref{con1} is equivalent
to a problem if any presymplectic manifold is contact.
So it might be a reasonable strategy to look for a presymplectic
form such that its Reeb foliation has properties of contact
forms. We consider this note as a step in this direction.

The following theorem gives a property which is shared by all
Reeb vector fields of contact forms. Before stating the theorem, let us recall
the notion of {\it basic cohomology.}

If  $M$ is a  closed manifold equipped with a foliation
$\mathfrak{F},$ then consider a cochain complex $(C^n,d^n)$ where
$C^n$ is the set of $n$-forms on M such that $$\iota_Y\alpha=0
\text{ and } \iota_Yd\alpha=0 $$ for any vector $Y$ tangent
to $\mathfrak{F}$ and $d^n$ is the usual exterior derivative. Then
by {\it n-th basic cohomology group} of $(M,\mathfrak F)$ we mean
the group  $H^n_b(M,\mathfrak{F})=\text{ker } d^n \slash \text{ im }
d^{n-1}.$

\begin{thm}\label{contcoh} Assume that $\lambda$ is a contact form on a closed
manifold  $M$ with its Reeb field $R.$ Then the basic cohomology
class $[d\lambda] \in H^2_b(M,R)$ is non-zero.
\end{thm}

{\bf Proof.} By the definition of Reeb field,
$\lambda (R) \equiv  1.$ Assume that $[d\lambda ]_b=0,$ so that there
exists $\alpha$ such that $\alpha (R)=0$ and $d\alpha = d\lambda .$
Thus  $\phi_0 = \alpha - \lambda$ is closed and equal to 1 on $R.$
It yields existence of a closed form $\phi$ which is $C^1$-close to
$\phi_0$ and such that $\ker \phi$ is integrable with compact
leaves by \cite{Ti}. On each leaf  $d\lambda$ restricts to a symplectic
form, so we get a contradiction, since the cohomology class of a symplectic form
on a closed manifold is non-zero.  \hfill\bk

Theorem \ref{contcoh} is no longer true for presymplectic forms. In
fact, one can have   $H^2_b(M,\mathfrak F)=0$ for a Reeb foliation
of a presymplectic form. Consider the following example studied by
Carri\`{e}re in \cite{C}.

\begin{ex}\label{car} {\em Let $T^3_A$ denotes the $T^2-$bundle over
the circle whose monodromy is given by matrix $A \in SL(2,\mathbb
Z)$ such that $tr A>2,$

$$T_A^3 = T^2 \times \mathbb R \slash (x,t) \sim (Ax,t+1).$$

\noindent Then both eigenvalues $\la, \frac{1}{\la}$ of $A$ are real and irrational.
Let $\mu_1, \mu_2$ be corresponding
eigenvectors. Define 1--forms $v_1,v_2$ on $T^2$ by
$v_1(\mu_1) = 1, \ v_1(\mu_2) = 0$ and $v_2(\mu_1) = 0, \ v_2(\mu_2) = 1.$
Extend $v_1,v_2$  to $T^2 \times \mathbb{R}$
by setting

$$\al_1= \la^t v_1$$
$$\al_2= \frac{1}{\la^t}v_2.$$ Now the forms $\al_1,\al_2$ are also well
defined on $T_A^3.$ By direct calculation, $d\al_1 = \ln(\la) dt
\wedge\al_1$ and $d\al_2 = -\ln(\la) dt\wedge\al_2.$ Thus  $d\al_1$
is a presymplectic form on $T^3_A$ with associated Reeb field $R =
\mu_2$ (and $d\alpha_2$ is presymplectic with the Reeb field $\mu_1 ).$
Furthermore, $\al_1(\mu_2)=0,$ hence by definition $[d\al_1]=0$
in $H^2_b(T_A^3,\mathfrak{R}).$

Carri\`{e}re also shows that the group
$H^2_b(T_A^3,\mathfrak{R})$ vanishes. \bk}
\end{ex}

Second property of vector fields we consider is related
to the condition $\eta (R)=1.$
Let  $\mathfrak V$ be a 1-dimensional foliation on $M.$

\begin{defin}\label{conformfol} A 1-form $\eta$ is a
{\it connection} on $\mathfrak V$ if
$\iota_Vd\eta=0$ and $\eta(V)\equiv 1$ for a vector field $V$ tangent
to $\mathfrak V.$
\end{defin}

The corresponding notion for vector fields is the following.

\begin{defin}\label{conformvec} Let $V$ be a nowhere vanishing vector field on
a manifold $M.$ We say that a 1-form $\eta$
is a connection form on $V$ if $\iota_Vd\eta=0$ and $\eta(V)\equiv 1.$
\end{defin}

Thus a 1-dimensional foliation has a connection if a vector field
tangent to the foliation has one. Note that any connection form
on $V$ is $V$-invariant, since $L_V\eta = \iota_Vd\eta + d\iota_V\eta = 0.$
If $V$ generates a (locally) free circle action, then what we have defined
becomes the standard notion of a connection of a principal bundle.

There are obstructions for a non-vanishing vector field $V$ to have a
connection form found by Dennis Sullivan \cite{S}.
Those are certain currents on $M$ that can be arbitrarily well approximated by
the boundary of a two-chain tangent to $V.$

\medskip
Vector fields described in Example \ref{car}
have no connection forms. To show this, we will
use the following observation.

\begin{lemma}\label{concon} If $ d\phi$ is a presymplectic form on
a manifold $M,$  $R$ is the Reeb field of $d\phi$ such
that $\phi(R)\geq 0$ and $\eta$
is a connection form for $R,$ then $\beta=K\phi+\eta$ is a contact form provided
that $K$ is large enough. Its Reeb vector field is $(K\phi (R)+1)^{-1}R.$
\end{lemma}

{\bf Proof.} To prove that $\beta$ is contact, it is enough to show
that $d(K\phi+\eta)= K d\phi+d\eta$
is non-degenerated on a subbundle $R^{\perp}$ transverse to $R$
and that $(K\phi+\eta)(R)>0.$ Since the second property is obvious,
we only show the first. Non-degeneracy is an open condition,thus
$d\phi +\frac 1K d\eta$ is non-degenerated if $K$ is
large enough. By assumptions, $\iota_R(Kd\phi+d\eta )=0,$
so $R$ is the Reeb field of $d\beta .$ Finally,
$\beta (R)\geq\eta (R)>0.$
 \bk

Lemma \ref{concon} implies that if the vector field $R$ described in Example \ref{car}
 had a connection form, then $R$ would be the Reeb field of a contact form
and by
Theorem \ref{contcoh} the basic cohomology class would be non-zero.
This is not possible, as $H^2_b(T_A^3,\mathfrak{R})$ is trivial.
 This contradicts also the Taubes theorem about
existence of closed orbits of a contact Reeb field (the Weinstein
conjecture) since trajectories of $R$ are irrational lines in $T^2.$

\medskip

Sullivan proved also that a
vector field $V$ admits a  connection if and only if it
is {\it geodesible,} i.e. there exists a
Riemannian metric
$g$ that makes the orbits of  $V$ geodesics. For  the proof
that every closed,
oriented odd-dimensional manifold has such field
we will use open book decompositions.

\begin{defin} An open book decomposition of $M$ is given by
\begin{enumerate}
\item a codimension two submanifold $B\subset M,$
\item a fibration $\pi:M \backslash B \rightarrow S^1$ with fibre $P,$
\item a tubular neighborhood $U$ of $B$ diffeomorphic to
$B \times D^2$
\end{enumerate}
such that the monodromy of the fibration $\pi$ is equal to the identity
in $P\cap U$ and   $\pi|U$ can be identified with the standard
projection $B \times (D^2\slash\{0\}) \rar S^1.$

\end{defin}

The submanifold $B$ is called the {\it binding} and the closure of $P$ is
the {\it page.}

The existence of an open book decomposition for an odd-dimensional closed
and oriented manifold was proved by Frank Quinn (\cite{Q}).
See also the discussion in Chapter 29 and in Appendix of
\cite{RW}.

\section{Existence of vector fields with Properties A,B}

By the preceding section,
any Reeb field $V$ of a contact form has Properties A and B.
In this section we prove that vector fields having both properties
exist on every closed oriented odd-dimensional
manifold.

\begin{thm}\label{exgeo} Any closed oriented odd-dimensional
manifold $M^{2n+1}$ admits a vector field with  connection
(hence a geodesible vector field).
\end{thm}

\NI The proof is based on the following lemma.

\begin{lemma}\label{mangeo} Consider an open book decomposition
of  a closed manifold $N.$ Then any vector field with connection on the binding
extends to a vector field with connection on $N.$

\end{lemma}

{\bf Proof.} By assumptions, we have a vector field $X_B$ on $B$ tangent to $B$ together with a connection form $\alpha$ for $X_B.$ The mapping $\pi:B \times (D^2\slash\{0\}) \rar S^1$ can be given in polar coordinates by $\pi(b,(r,\varphi))=\varphi,$ hence we set a connection form
$\pi^*d\varphi$
for the vector field $\frac{\partial}{\partial \varphi}$ on the boundary of
$B \times D^2.$
This pair can easily be extended to $N \backslash (B \times D^2)$
by taking any lift of $\frac{\partial}{\partial \varphi}$ that coincides on
the boundary $\partial(B \times D^2)$ with $\frac{\partial}{\partial \varphi},$ while the connection form for $\frac{\partial}{\partial \varphi}$ on $N \backslash (B \times D^2)$ remains the same and is equal to $\pi^*d\varphi$ . The proof is completed by showing that for suitably chosen functions $f,g:[0;1]\rightarrow [0;1]$
($f,g=0$ for $r \in [0;\varepsilon)$ and $f,g=1$ for $r \in (1-\varepsilon;1]$)
the form
$$\eta=f(r)d\varphi+(1-f(r))\alpha$$
is a connection form for
$$V=g(r)\frac{\partial}{\partial \varphi}+(1-g(r))X_B.$$
We have
$$\eta(X)=f(r)g(r)+(1-f(r))(1-g(r))>0,$$
and furthermore
$$d\eta = f'(r)drd\varphi-f'(r)dr\wedge\alpha+(1-f(r))d\alpha,$$
hence
$$\iota_Vd\eta=(-f'(r)g(r)+f'(r)(1-g(r)))dr=f'(r)(1-2g(r))dr.$$
This means that if $f$ is not constant, then $g(r)=\frac 12.$ Such functions
are easy to construct.
\bk

{\bf Proof} of Theorem \ref{exgeo}. The proof is by induction on $n.$ For $n=1$ this
follows from the fact that every closed oriented 3-manifold is contact
(c.f. \cite{TW},\cite{Ma}), and it is trivial for $n=0.$
In general, M has an open book decomposition (\cite{Q}) and we apply Lemma \ref{mangeo}
to complete the proof. \bk

Now we deal with Property B. The following lemma will be
useful for recognizing whether, for a connection form $\eta$ on a vector field
$V,$  the class
$[d\eta]_b \in H^2_b(M,\mathfrak{V})$  is non-zero.

\begin{lemma}\label{basiccoh} Let $\eta$ be a connection form for a vector field $V.$
If $V$ has a closed orbit $\gamma$ which is
homologically trivial and $\int_{\gamma}\eta \neq 0,$ then
$[d\eta]_b \neq 0.$ The same is true if there are two closed
orbits $\gamma , \gamma'$ of $V$ field representing the same
homology class and $\int_{\gamma}\eta \neq \int_{\gamma'}\eta .$
\end{lemma}

{\bf Proof.} If $[d\eta]_b=0,$  then we have 1-form $\phi$ such
that $d\phi=d\eta$ and $\phi(\mathcal R)=0.$ Consequently,
$\eta-\phi$ is closed and $\int_{\gamma}(\eta-\phi)\neq 0,$
contrary to homological triviality of $\gamma.$ This argument works also
for two closed orbits with different values of $\int_{\gamma}\eta
.$ \bk

\begin{cor}\label{main} Every closed oriented odd-dimensional manifold $M^{2n-1}$
has a geodesible vector field $V$ with such a connection form $\eta$
that $[d\eta]_b \neq 0$ in $H^2_b(M,\mathfrak{V}).$
\end{cor}

{\bf Proof.} By Theorem \ref{exgeo}, it is sufficient to prove
that for a vector field $V$ constructed in Theorem \ref{exgeo} we can
choose $\eta$ such that  $[d\eta]_b \neq 0.$ Observe that near the boundary
$\partial (B \times D^2)$ we have $V = \frac{\partial}{\partial t},$
hence it has a contractible closed
orbit $\gamma.$ The connection form $\eta$ given by Lemma \ref{mangeo} in this subset
is $dt,$ hence $\int_{\gamma}\eta \neq 0.$ Lemma \ref{basiccoh} completes the proof.
\bk

\begin{rem} {\em For presymplectic manifolds $M^{2n+1}$ Theorem \ref{main}
can be derived directly from \cite{MMP}, where the authors proved that $M$
has an open book decomposition with presymplectic binding. As in Theorem \ref{main},
we proceed with induction on $n.$}
\end{rem}

\section{Presymplectic confoliation}

Assume that a codimension one subbundle $\xi$ on $M^{2n+1}$ is equal to
the kernel of a 1-form $\alpha .$  If
$\alpha \wedge (d\alpha)^n \geq 0,$ then $\xi$ is called a (positive)
{\it confoliation} on $M.$ The form is determined by $\xi$ up to a positive function
and is contact if $\alpha \wedge (d\alpha)^n > 0.$
A form defining a confoliation will be called confoliation form.

\begin{ex}\label{conpre} {\em Let $S^1\rar M^{2n+1}\rar B$ be a principal $S^1$
fibration and
$\eta$ an invariant connection
such that its curvature defines a symplectic form on $B.$  Then $\eta$ is
the Boothby-Wang contact form \cite{BW}. If we assume only that
$(d\eta )^n\geq 0,$ then we get a confoliation form.

Of course, any closed nowhere vanishing 1-form $\alpha$ is confoliation form,
but then $\xi$ is integrable and if $M$ is closed, then it fibers over a circle.
More exactly, there is a nowhere vanishing 1-form $\alpha'$ which is closed,
$C^1$-close to $\alpha$  and $\xi' = \ker \alpha'$ is tangent
to fibers of the fibration over $S^1$ (\cite{Ti}).} \bk
\end{ex}

Some methods to deform a confoliation form to a contact form were
found,
see for example second chapter of \cite{ET} for dimension 3 case and
\cite{AW} for the general case. Generally,
there are two types of obstacles
for a confoliation form  $\alpha$ to be contact. First of all the
rank of $d\alpha$ can be strictly less then $2n.$ Even if the rank
of $d\alpha$ is maximal (hence $d\alpha $ is presymplectic),
then its Reeb field $R$ might
lie in $\ker \alpha$ (thus $\alpha\wedge (d\alpha )^n=0)$ at some points.

\begin{defin} A 1-form $\alpha$ is
a presymplectic confoliation form
if $\alpha$ is a (positive) confoliation form and $d\alpha$ is presymplectic.
\end{defin}

Direct calculation shows that conditions defining presymplectic
confoliation form are equivalent to the following properties of
$\alpha:$ $d\alpha$ is presymplectic and $\alpha(R)\geq 0.$
Under the assumption that  the Reeb field admits a connection,
one can deform a presymplectic confoliation form into a contact form by
Lemma \ref{concon}.

\begin{prop} Assume that $\alpha$ is a presymplectic confoliation
form and the Reeb field of $d\alpha$ has a connection form $\eta .$
Then for $\varepsilon$ small enough the
form $\alpha+\varepsilon\eta$ is contact.
\end{prop}

\begin{ex} {\em The forms $\alpha_1, \alpha_2$ defined  in Example \ref{car}
are presymplectic confoliation forms, but
their Reeb fields have no connections.
However, the form
$\phi_{\varepsilon} = \al_1 + \varepsilon \al_2$ is contact for
$\varepsilon >0$ since
$\phi_{\varepsilon} \wedge d\phi_{\varepsilon} = 2\varepsilon
\ln(\la) \al_1 \wedge\al_2 \wedge dt > 0.$ Recall that exactly one of
eigenvalues $\lambda,\frac{1}{\lambda}$ of $A$ is greater than 1,
thus we must choose suitable $\alpha_i$ to have that $\phi_{\varepsilon}
\wedge d\phi_{\varepsilon}$ is positive; therefore we take $\alpha_1$
from Example \ref{car} and perturb it linearly in the direction of $\alpha_2.$

Furthermore, contact forms $\phi_{\varepsilon}$ are $C^1$-approximate
of the presymplectic confoliation form
 $\al_1.$} \bk
\end{ex}

\bibliographystyle{amsalpha}

\medskip

\noindent {\bf Mathematical Institute, Wroc\l aw University,

\noindent pl. Grunwaldzki 2/4,

\noindent 50-384 Wroc\l aw, Poland}

\medskip

\NI and

\medskip

\noindent {\bf Department of Mathematics and Information Technology,

\noindent University of Warmia and Mazury,

\noindent \.{Z}o\l nierska 14A, 10-561 Olsztyn, Poland}

\medskip

\begin{flushleft}
\tt hajduk@math.uni.wroc.pl
\end{flushleft}
\medskip

\noindent {\bf West Pomeranian University of Technology,

\noindent Mathematical Institute}

\noindent {\bf Al. Piast\'{o}w 48/49, 70--311 Szczecin, Poland}
\medskip

\begin{flushleft}
 \tt rafal\_walczak2@wp.pl
\end{flushleft}

\end{document}